\documentclass[12pt,a4paper]{article}
\usepackage{pifont}

\usepackage{epsfig}
\usepackage{amsmath}
\usepackage{amssymb}
\usepackage{graphicx}

\newtheorem{propo}{Proposition}[section]

\newtheorem{lemma}[propo]{Lemma}

\newtheorem{theo}[propo]{Theorem}
\newtheorem{examp}[propo]{Example}

\newtheorem{prop}[propo]{Proposition}
\newtheorem{rem}[propo]{Remark}
\newtheorem{constr}[propo]{Construction}

\newcommand{\bl}{\begin{lemma}}
\newcommand{\el}{\end{lemma}}

\def\d12{{_{12}}}

 \def\a{\alpha} \def\b{\beta} 
\def\a{\alpha} \def\b{\beta}

\def\A{{\rm Aut}}

\def\Cay{{\rm Cay}}

\def\Cos{{\rm Cos}}

\usepackage{indentfirst,latexsym,bm}

\def\A{\rm A}

\def\Aut{{\rm Aut}}

\def\a{\alpha}
\def\b{\beta}

\def\diam{{\rm diam}}

\def\K{{\rm K}}

\def\H{{\rm H}}
\def\girth{{\rm girth}}

\def\geod{{\rm geod}}

\begin{document}
\title{On distance, geodesic and arc transitivity of graphs
}

\author{Alice Devillers, Wei Jin\thanks{The  second author is  supported by
the Scholarships for International Research Fees (SIRF) at UWA.},
Cai Heng Li and Cheryl E. Praeger\thanks{ This paper  forms part of
Australian Research Council Federation Fellowship FF0776186 held by
the fourth author. The first author is supported by UWA as part of
the Federation Fellowship project. }\thanks{E-mail addresses:
alice.devillers@uwa.edu.au (A.Devillers), weijin@maths.uwa.edu.au
(W.Jin), cai.heng.li@uwa.edu.au (C.H.Li) and
cheryl.praeger@uwa.edu.au (C.E.Praeger).  }
\\{\footnotesize
Centre for the Mathematics of Symmetry and Computation, School of
Mathematics and Statistics,  }
\\{\footnotesize
The University of Western Australia, Crawley, WA 6009, Australia} }

\date{ }

\maketitle

\begin{abstract}

We compare three transitivity properties of finite graphs, namely,
for a positive integer $s$, $s$-distance transitivity, $s$-geodesic
transitivity and $s$-arc transitivity. It is known that if a finite
graph is $s$-arc transitive but not $(s+1)$-arc transitive then
$s\leq 7$ and $s\neq 6$. We show that there are infinitely many
geodesic transitive graphs with this property for each of these
values of $s$, and that  these graphs can have arbitrarily large
diameter if and only if $1\leq s\leq 3$. Moreover, for a prime $p$
we prove that there exists a graph of valency $p$ that is 2-geodesic
transitive but not 2-arc transitive if and only if $p\equiv  1\pmod
4$, and for each  such prime there is a unique graph with this
property: it is an antipodal double cover of the complete graph
$K_{p+1}$ and is geodesic transitive with automorphism group
$PSL(2,p)\times Z_2$.

\end{abstract}

\vspace{2mm}

 \hspace{-17pt}{\bf Keywords:} Graphs; distance-transitivity; geodesic-transitivity; arc-transitivity

\section{Introduction}

The  study of finite distance-transitive graphs goes back  to
Higman's paper \cite{DGH-1} in which ``groups of maximal diameter"
were introduced. These are permutation groups which act distance
transitively on some graph. The family of distance transitive graphs
includes many interesting and important graphs, such as the Johnson
graphs, Hamming graphs, Odd graphs, Paley graphs and certain
antipodal double covers of complete graphs that are discussed in
this paper.

We examine graphs with various symmetry properties   which are
stronger than arc-transitivity. The weakest of these properties is
$s$-distance transitivity where $s$ is at most the diameter of the
graph (see Section 2 for a precise definition), in which,  all pairs
of vertices at a given distance at most $s$ are equivalent under
graph automorphisms. If $s$ is equal to the diameter, the graph is
said to be \emph{distance transitive}.

For a fixed integer   $s\geq 2$, two stronger concepts than
$s$-distance transitivity  are important for our work. The first is
$s$-geodesic transitivity for $s$ at most the diameter, in which for
each integer $t\leq s$, all ordered $t$-paths $(v_0,v_1,\cdots,v_t)$
for which $v_0,v_t$ are at distance $t$, are equivalent under graph
automorphisms. Such $t$-paths are called \emph{$t$-geodesics}. The
second property is $s$-arc transitivity, in which for each integer
$t\leq s$, all ordered $t$-walks $(v_0,v_1,\cdots,v_t)$, with
$v_{i-1}\neq v_{i+1}$, for each $i=1,2,\cdots, t-1$, are equivalent.
It is straightforward to verify that, for $s$ at most the diameter,
$s$-arc transitivity implies $s$-geodesic transitivity, which in
turn implies $s$-distance transitivity. The purpose of this paper is
to provide some insight into the differences between these
conditions, especially for $s=2$.

A  graph is called  \emph{geodesic transitive} if it is $s$-geodesic
transitive for $s$ equal to the diameter. Geodesic transitive graphs
are in particular  1-arc transitive, and may or may not be $s$-arc
transitive for some $s>1$. Our first result specifies the possible
arc transitivities for geodesic transitive graphs. We note that (see
\cite{weiss}) if a finite graph is $s$-arc transitive but not
$(s+1)$-arc transitive then $s\in \{1,2,3,4,5,7\}$.

\begin{theo}\label{gtnotat}
For each $s\in \{1,2,3,4,5,7\}$, there are infinitely many geodesic
transitive graphs that are $s$-arc transitive but not $(s+1)$-arc
transitive. Moreover, there exist geodesic transitive graphs that
are $s$-arc transitive but not $(s+1)$-arc transitive with
arbitrarily large diameter if and only if $s\in \{1,2,3\}$.
\end{theo}

Theorem \ref{gtnotat} is proved in Section 3 by analysing some
well-known families of distance transitive graphs, namely, Johnson
graphs, Hamming graphs, Odd graphs and classical generalized
polygons. In fact if $s\geq 4$, then all graphs with the property of
Theorem \ref{gtnotat} have diameter at most 8 and are known
explicitly, see Proposition \ref{dt-4at-diam}.


\medskip

Our second result is the main result of the paper, proved in
Subsection 5.1. It classifies explicitly all $2$-geodesic transitive
graphs of prime valency  that are not 2-arc transitive.

\medskip

\begin{theo}\label{gt-primeval}
Let   $\Gamma$ be a connected $2$-geodesic transitive but not
$2$-arc transitive graph of prime valency $p$. Then $\Gamma$ is a
nonbipartite antipodal  double cover of the complete graph
$K_{p+1}$, where $p\equiv  1\pmod 4$. Further, $\Gamma$ is geodesic
transitive and $\Gamma$ is unique up to isomorphism with
automorphism group $PSL(2,p)\times Z_2$.
\end{theo}

This family of graphs arose also  in
\cite{GLP,Li-ci-soluble,Taylor-1} in different contexts, and they
are distance transitive of diameter 3. We prove in Lemma
\ref{val-p-lemma-4} that each graph in the family  is geodesic
transitive. In Subsection 5.2, we construct these graphs as in
\cite{Li-ci-soluble} and prove that they satisfy the hypotheses of
Theorem \ref{gt-primeval}. It would be interesting to know if a
similar classification is possible for non-prime valencies. This is
the subject of further research by the second author.

We  complete our comparison of these transitivity properties by
producing 2-distance transitive graphs that are not 2-geodesic
transitive. We give just one infinite family of examples, namely the
Paley graphs $P(q)$ where $q\geq 13$ is a prime power and $q\equiv 1
\pmod{4}$ (see Section 4 for a definition). These graphs $P(q)$ have
diameter 2 and are well-known to be distance transitive. The
information about these graphs is important for the proof of Theorem
\ref{gtnotat}.

\begin{theo}\label{exam-valp-ha}
Let  $q\equiv 1 \pmod{4}$ be a prime power. Then the Paley graph
$P(q)$ is distance transitive for all $q$, but $P(q)$ is geodesic
transitive if and only if $q=5$ or $9$.
\end{theo}

\begin{rem}\label{dnotg-rem-1}
{\rm The Paley graphs $P(q)$ with $q>9$ seem to be  the first family
of diameter $2$ graphs observed to be distance transitive but not
geodesic transitive. Since all diameter $2$ distance transitive
graphs are known, it would be interesting to know which of them are
geodesic transitive.}

\end{rem}

These results give some insight into the relationship between
$s$-distance transitivity, $s$-geodesic transitivity and $s$-arc
transitivity for $s=2$. It would be interesting to understand the
relationship between these properties for larger values of $s$.

\section{Preliminaries}

All graphs of this paper are finite, undirected simple graphs. We
first give some definitions which will be used throughout the paper.
Let $\Gamma$ be a graph. We use $V\Gamma,E\Gamma$, and $\Aut\Gamma$
to denote its \emph{vertex set}, \emph{edge set} and \emph{ full
automorphism group}, respectively. The size of $V\Gamma$ is called
the \emph{order} of the graph.  The graph $\Gamma$ is said to be
\emph{vertex transitive} (or \emph{edge transitive}) if the action
of $\Aut \Gamma$ on $V\Gamma$ (or $E\Gamma$) is transitive.

For (not necessarily distinct) vertices $u$ and $v$ in $V\Gamma$, a
$walk$ from $u$ to $v$ is a finite sequence of vertices
$(v_{0},v_{1},\cdots,v_{n})$ such that $v_{0}=u$, $v_{n}=v$ and
$(v_{i},v_{i+1})\in E\Gamma$ for all $i$ with $0\leq i<n$, and $n$
is called the \emph{length} of the walk. If  $v_{i}\neq v_{j}$ for
$0\leq i < j\leq n$, the  walk is called a \emph{path} from $u$ to
$v$. The smallest value for $n$ such that there is a path of length
$n$ from $u$ to $v$ is called the \emph{distance} from $u$ to $v$
and is denoted $d_{\Gamma}(u, v)$.  The \emph{diameter}
$\diam(\Gamma)$ of a connected graph $\Gamma$ is the maximum of
$d_{\Gamma}(u, v)$ over all $u, v \in V\Gamma$.

Let $G \leq \Aut\Gamma$ and $s\leq \diam(\Gamma)$. We say that
$\Gamma$ is \emph{$(G, s)$-distance transitive} if, for any two
pairs of vertices $(u_1,v_1)$, $(u_2,v_2)$ with the same distance
$t\leq s$,  there exists $g\in G$ such that $(u_1,v_1)^g=(u_2,v_2)$.

For a positive integer $s$, an $s$-arc of $\Gamma$ is a walk
$(v_0,v_1,\cdots,v_s)$ of length $s$ in  $\Gamma$ such that
$v_{j-1}\neq v_{j+1}$ for $1\leq j\leq s-1$. Moreover, a 1-arc is
called an \emph{arc}. Suppose $G\leq \Aut \Gamma$. Then $\Gamma$ is
said to be \emph{ $(G,s)$-arc transitive} if, for any two $t$-arcs
$\alpha$ and $\beta$ where $t\leq s$,  there exists $g\in G$ such
that $\alpha^g=\beta$.  A remarkable result of Tutte about
$(G,s)$-arc transitive graphs with valency three shows that $s \leq
5$, see \cite{Tutte-1,Tutte-2}. About twenty years later, relying on
the classification of  finite simple groups, Weiss in \cite{weiss}
proved that there are no  $(G,8)$-arc transitive graphs with valency
at least three.   For more work on $(G,s)$-arc transitive graphs see
\cite{IP-1,Praeger-4,Praeger-1993-1}.




Let $u,v$ be  distinct vertices of $\Gamma$. Then a path of shortest
length from $u$ to $v$ is called a \emph{geodesic} from $u$ to $v$,
or sometimes  an \emph{$i$-geodesic} if $d_{\Gamma}(u,v)=i$.
Moreover, 1-geodesics are  arcs.  If $\Gamma$ is connected, for each
$i \in \{1,\cdots,\diam(\Gamma) \}$, we set $geod_i(\Gamma)=\{$all
$i$-geodesics  of $\Gamma\}$. Let $G\leq \Aut\Gamma$  and $ s \leq
\diam(\Gamma)$. Then $\Gamma$ is said to be \emph{ $(G,s)$-geodesic
transitive} if,  for each $i=1,2,\cdots,s$, $G$ is transitive on
$\geod_i(\Gamma)$. When $s=\diam(\Gamma)$, $\Gamma$ is said to be
\emph{$G$-geodesic transitive}. Moreover,  if we do not wish to
specify the group we will say that $\Gamma$ is \emph{ $s$-geodesic
transitive} or \emph{ geodesic transitive} respectively, and
similarly for the other properties.

The following are some examples of geodesic transitive graphs.

\begin{examp} \label{geod-deg-exam}
{\rm(i)} For any $n\geq 1$, both the complete graph $\K_n$ and the
complete bipartite graph $K_{n,n}$ are geodesic transitive.

{\rm (ii)} Let $\Gamma=\K_{m[b]}$ be a complete multipartite graph
with $m\geq 3$ parts of size $b\geq 2$. Then $A:=\Aut \Gamma=S_b\wr
S_m$ is transitive on $V\Gamma$ and on the set of arcs
$\geod_1(\Gamma)$. Let $(u,v)$ be an arc of $\Gamma$. Then
$|\Gamma_2(u)\cap \Gamma(v)|=b-1$ and  $A_{u,v}$ induces $S_{b-1}$
on $\Gamma_2(u)\cap \Gamma(v)$, so $\Gamma$ is $2$-geodesic
transitive. Since the diameter of $\Gamma$ is $2$, it follows that
$\Gamma$ is geodesic transitive.
\end{examp}


We also give three infinite families of geodesic transitive graphs
with arbitrarily large diameter in Section 3.

If  a graph $\Gamma$ is $(G,s)$-arc transitive and  $s\leq
\diam(\Gamma)$, then  $s$-geodesics and  $s$-arcs are same, and
$\Gamma$ is $(G,s)$-geodesic transitive. However, $\Gamma$ can be
$(G,s)$-geodesic transitive but not $(G,s)$-arc transitive. The
girth of $\Gamma$, denoted by $girth(\Gamma)$, is the length of the
shortest cycle in $\Gamma$.
\begin{lemma}\label{agd-lemma-1}
Suppose that a graph $\Gamma$ is $(G,s)$-geodesic transitive for
some $G\leq \Aut \Gamma$ with $2\leq s\leq \diam(\Gamma)$. Then
$\Gamma$ is $(G,s)$-arc transitive if and only if
$\girth(\Gamma)\geq 2s$.

\end{lemma}
{\bf Proof.} Note that each $i$-geodesic is an $i$-arc, for  $1\leq
i\leq \diam(\Gamma)$. Thus $\Gamma$  is $(G,s)$-arc transitive if
and only if  each $s$-arc is an $s$-geodesic, and this is true  if
and only if $\girth(\Gamma)\geq 2s$. $\Box$

\bigskip

An example of graphs which do not have the property of Lemma
\ref{agd-lemma-1} for $s=2$ are  the complete multipartite graphs
$\K_{m[b]}$ with $m\geq 3$ parts of size $ b\geq 2$, which have
girth 3,  are 2-geodesic transitive (see Example \ref{geod-deg-exam}
(ii)), but are  not 2-arc transitive (by Lemma \ref{agd-lemma-1}).


\medskip

In our study of Paley graphs we use the concept of a Cayley graph.
For a finite group $G$, and a subset $S$ of $G$ such that $1\notin
S$ and $S=S^{-1}$, the \emph{Cayley graph} $\Cay(G,S)$ of $G$ with
respect to $S$ is  the graph with vertex set $G$ and edge set
$\{\{g,sg\} \,|\,g\in G,s\in S\}$. The group $R(G)=\{\rho_x|x\in
G\}$, where $\rho_x:g\mapsto gx$,  is a subgroup of the automorphism
group of $\Cay(G,S)$ and acts regularly on the vertex set, that is
to say, $R(G)$ is transitive and only the identity $\rho_{1_G}$
fixes a vertex. It follows that $\Cay(G,S)$ is vertex transitive.

The following is a criterion for a connected graph to be a Cayley
graph.
\begin{lemma}{\rm(\cite[Lemma 16.3]{Biggs-1})}\label{cayley-1}
Let $\Gamma$ be a connected graph. Then a subgroup $H$ of $\Aut
\Gamma$ acts regularly on the vertices if and only if $\Gamma$ is
isomorphic to a Cayley graph $\Cay(H,S)$ for some set $S$ which
generates $H$.
\end{lemma}

For a graph $\Gamma$, the \emph{k-distance graph} $\Gamma_k$ of
$\Gamma$ is the graph with vertex set $V\Gamma$, such that  two
vertices are adjacent if and only if they are at distance $k$ in
$\Gamma$. If $d=\diam(\Gamma)\geq 2$, and $\Gamma_d$ is a disjoint
union of complete graphs, then  $\Gamma$ is said to be an
\emph{antipodal graph}.

Suppose that $\Gamma$ is an antipodal distance-transitive graph of
diameter $d$. Then we may partition its vertices into sets, called
\emph{fibres}, such that any two distinct vertices in the same fibre
are at distance $d$ and two vertices in different fibres are at
distance less than $d$. Godsil, Liebler and Praeger gave a complete
classification of  antipodal distance transitive covers of complete
graphs.

The following lemma follows directly from the Main Theorem of
\cite{GLP}. We shall apply it to characterise the examples in
Theorem \ref{gt-primeval}.

\begin{lemma}{\rm(\cite{GLP})}\label{anticover-lemma-1}
Suppose that $G$ is a distance transitive automorphism group of a
finite nonbipartite graph  $\Gamma$. Suppose further that $\Gamma$
is antipodal with fibres of size $ 2$ and with antipodal quotient
the complete graph $K_n$.  Then either $\Gamma=K_{n[2]}$ of diameter
$2$, or $\Gamma$ has diameter $3$ and one of the following holds.

{\rm I.} $\Gamma$ appears in {\rm\cite{Taylor-1}} and is one of

\ \ \ \ {\rm (a)} $n=2^{2m-1}\pm 2^{m-1}, G\leq 2\times Sp(2m,2)$
for some $m\geq 3$.

\ \ \ \ {\rm (b)} $n=2^{2a+1}+1$, $G\leq 2\times \Aut(R(q))$ for
some $a\geq 1$.

\ \ \ \ {\rm (c)} $n=176$, $HiS\leq G\leq 2\times HiS$, or $n=276$,
$Co_3\leq G\leq 2\times Co_3$.

\ \ \ \ {\rm (d)} $n=q^3+1$, $G\leq 2\times P\Gamma U(3,q^2)$ for
some  $q>3$.

\ \ \ \ {\rm (e)} $n=q+1$, $G\leq 2\times P\Sigma L(2,q)$ for some
$q\equiv 1 \pmod 4$.

{\rm II.} $\Gamma$ is a graph appearing in Example {\rm 3.6} of
{\rm\cite{GLP}},  and $n=q^{2d}$ with $q$ even.

\end{lemma}

\medskip

\noindent \textbf{Partitions and quotient graphs:} \quad  Let $G$ be
a group of permutations acting on a set $\Omega$. A
\emph{$G$-invariant partition} of $\Omega$  is a partition
$\mathcal{B}=\{B_1,B_2,\cdots,B_n\}$ such that for each $g\in G$,
and each $B_i \in \mathcal{B}$, the image  $B_{i}^g\in \mathcal{B}$.
The parts of $\Omega$ are often called \emph{blocks} of $G$ on
$\Omega$. For a $G$-invariant partition $\mathcal{B}$ of $\Omega$,
we have two smaller transitive permutation groups, namely the  group
$G^{\mathcal{B}}$ of permutations of $\mathcal{B}$ induced by $G$;
and the group $G_{B_i}^{B_i}$ induced on $B_i$ by $G_{B_i}$ where
$B_i\in \mathcal{B}$.

Let $\Gamma$ be a graph, and $G\leq \Aut \Gamma$. Suppose
$\mathcal{B}=\{B_1,B_2,\cdots,B_n \}$ is a $G$-invariant partition
of $V\Gamma$. The \emph{quotient graph} $\Gamma_{\mathcal{B}}$ of
$\Gamma$ relative to $\mathcal{B}$ is defined to be the graph with
vertex set $\mathcal{B}$ such that $\{B_i,B_j\}$ is an edge of
$\Gamma_{\mathcal{B}}$ if and only if there exist $x\in B_i, y\in
B_j$ such that $\{x,y\} \in E\Gamma$. We say that
$\Gamma_{\mathcal{B}}$ is \emph{nontrivial} if $1< |\mathcal{B}|<
|V\Gamma|$.  The graph $\Gamma$ is said to be a \emph{cover} of
$\Gamma_{\mathcal{B}}$ if for each edge $\{B_i,B_j\}$ of
$\Gamma_{\mathcal{B}}$ and $v\in B_i$, we have $|\Gamma(v)\cap
B_j|=1$.

\section{Proof of Theorem \ref{gtnotat}}

In this section, we first describe three families of geodesic
transitive graphs, each with unbounded diameter and valency, namely
the Johnson graphs, Hamming graphs and Odd graphs. Graphs in these
families are $s$-arc transitive but not $(s+1)$-arc transitive for
various $s\leq 3$. In the last subsection, we  give the proof of
Theorem \ref{gtnotat}.

In the following  discussion, for  integers $i,j$, we define
$[i,j]=\{n\in Z\,|\,i\leq n\leq j\}$. Note that $[i,i]=\{i\}$ and
$[i,j]=\emptyset$ when $i>j$.

\subsection{Johnson graphs}

Let $\Omega=[1, n]$ where $n\geq 3$, and let $1\leq k\leq
[\frac{n}{2}]$ where $[\frac{n}{2}]$ is the integer part of
$\frac{n}{2}$.  Then the \emph{Johnson graph} $J(n,k)$ is the graph
whose vertex set is the set of all $k$-subsets of $\Omega$, and two
vertices $u$ and $v$ are adjacent if and only if $|u\cap v|=k-1$.
Let $\Gamma=J(n,k)$. By \cite[Section 9.1]{BCN},  $\Gamma$ has the
following properties: $\Gamma$ has diameter $k$, valency $k(n-k)$,
$\Aut \Gamma\cong S_n\times Z_2$ when $n=2k\geq 4$, otherwise $\Aut
\Gamma \cong S_n$, $\Gamma$ is distance transitive, and for any two
vertices $u$ and $v$,

\medskip
$u\in \Gamma_j(v)$ where $j\leq k$ if and only if $|u\cap v|=k-j$. \
\, \, \ \ \  $(J*)$

\medskip
Note that for $k=1$, $J(n,k)\cong K_n$ which has diameter 1. So in
the following discussion, we assume that $k\geq 2$ and $n\geq 4$.

\begin{prop}\label{j(n,k)}
Let $\Gamma=J(n,k)$ where $2\leq k\leq [\frac{n}{2}]$ and $n\geq 4$.
Then $\Gamma$ has girth $3$,  is geodesic transitive, but not
$2$-arc transitive.
\end{prop}
{\bf Proof.}  Since  $k\geq 2$ and $n\geq 4$, it follows that
$u_1=[1,k]$, $u_2=\{1\}\cup [3,k+1]$ and $u_3=[2,k+1]$ are three
vertices of $\Gamma$. By $(J*)$, they  are pairwise adjacent, so
$\Gamma$ has girth 3. Hence $\Gamma$ is not 2-arc transitive. We
will prove that $\Gamma$ is geodesic transitive. Since  $\Gamma$ is
distance transitive, it follows that $\Gamma$ is 1-geodesic
transitive.

Now suppose  that  $\Gamma$ is $(j-1)$-geodesic transitive where
$j\in [2,k]$. Let $\mathcal{V} =(v_0,v_1,\cdots,v_{j-1},v_j)$ where
$v_i=[1,k-i] \cup [k+1,k+i]$ for each $i\in [0,j]$. Then by $(J*)$,
$\mathcal{V}$ is a $j$-geodesic. Let $\mathcal{U}$ be any other
$j$-geodesic. Then since $\Gamma$ is $(j-1)$-geodesic transitive,
there exists $\a \in \Aut \Gamma$ such that
$\mathcal{U}^\a=(v_0,v_1,\cdots,v_{j-1},u_j)$ which is also a
$j$-geodesic, and in particular $u_j\in \Gamma_j(v_0)$, so $|v_0\cap
u_j|=k-j$ by $(J*)$. Since $v_{j-1}$ and $u_j$ are adjacent,
$|v_{j-1}\cap u_j|=k-1$. Hence there exist a unique $x$ such that
$\{x\}=v_{j-1}\setminus u_j$ and a unique $y$ such that
$\{y\}=u_{j}\setminus v_{j-1}$. First, if $x\geq k+1$, then
$[1,k-(j-1)]\subseteq u_j$, so $|v_0\cap u_j|\geq k-(j-1)$ which
contradicts  $|v_0\cap u_j|= k-j$. Thus $x\in [1,k-(j-1)]$, and
hence $[k+1,k+(j-1)]\subseteq u_j$. Second, if $y\leq k$, then $y\in
[k-(j-2),k]$. It follows that $v_0\cap u_j$ contains
$([1,k-(j-1)]\cup \{y\})\setminus \{x\}$, a set of size $k-(j-1)$,
which also contradicts $|v_0\cap u_j|= k-j$. Thus $y
> k$, and hence $y\in [k+j,n]$. Therefore,
$u_j=([1,k-(j-1)]\setminus \{x\})\cup [k+1,k+(j-1)]
 \cup \{y\}$. Let $A=\Aut \Gamma$.  Since
$[1,k-(j-1)]\subseteq v_0\cap v_1\cap \dots \cap v_{j-1}$ and
$[k+j,n]\subseteq \Omega \setminus (v_0\cup v_1\cup \dots \cup
v_{j-1})$, and since $Sym(\Omega)\leq A $, it follows that $Sym([1,
k-(j-1)])\times Sym([k+j,n]) \leq A_{v_0,v_1,\cdots,v_{j-1}}$. Hence
there exists $\b \in A_{v_0,v_1,\cdots,v_{j-1}}$ such that
$x^\b=k-(j-1)$ and $y^\b=k+j$, and so
$(v_0,v_1,v_2,\cdots,u_j)^\b=(v_0,v_1,v_2,\cdots,v_{j-1},v_j)=\mathcal{V}$,
that is  $\mathcal{U}^{\a \b}=\mathcal{V}$. This completes the
induction. Thus $\Gamma$ is geodesic transitive. $\Box$

\subsection{Hamming graphs}

Let  $n\geq 2$ and let $d$ be a positive integer. Then the
\emph{Hamming graph} $\H(d,n)$  has vertex set $Z_n^d=Z_n\times
Z_n\times \cdots \times Z_n$, seen as a module on the ring
$Z_n=[0,n-1]$, and two vertices $u,v$ are adjacent if and only if
$u-v$ has exactly one non-zero entry. For a vertex $u\in V\H(d,n)$,
we denote by $|u|$ the number of its non-zero entries. Then by
\cite[Section 9.2]{BCN}, $\H(d,n)$ has diameter $d$, valency
$d(n-1)$, is distance transitive, $\Aut \H(d,n)\cong S_n\wr S_d$,
and for two vertices $u,v$,

\medskip

$u\in \Gamma_k(v)$ where $k\leq d$ if and only if $|u-v|=k$. \ \, \,
\ \ \  $(H*)$

\begin{prop}\label{Ham-2-trans}
Let $\Gamma=\H(d,n)$ with  $d\geq 2$, $n\geq 2$. Then $\Gamma$ is
geodesic transitive. Moreover, if $n=2$ and $d\geq 3$ then $\Gamma$
has girth $4$ and is $2$-arc transitive but not $3$-arc transitive,
while if $n\geq 3$ then $\Gamma$ has girth $3$ and is  not $2$-arc
transitive.
\end{prop}
{\bf Proof.} Since $\Gamma$ is distance transitive,  it follows that
$\Gamma$ is 1-geodesic transitive.

Now, suppose  that  $\Gamma$ is $(j-1)$-geodesic transitive where
$j\in [2,d]$.  Let $\mathcal{V}=(v_0,v_1,\cdots,v_j)$ where for each
$i\in [0,j]$, $v_i=(1,\cdots,1,0,\cdots,0)$,  the first $i$ entries
are equal to $1$ and the last $(d-i)$ entries are equal to 0. Then
by $(H*)$, $\mathcal{V}$ is a $j$-geodesic.

Suppose that $\mathcal{U}$ is any other $j$-geodesic of $\Gamma$.
Since  $\Gamma$ is $(j-1)$-geodesic transitive, there exists $\a \in
\Aut \Gamma$ such that $\mathcal{U}^\a=(v_0,v_1,\cdots,v_{j-1},u_j)$
for some $u_j$. Suppose that the last $d-(j-1)$ entries of $u_j$ are
0. Since $v_{j-1},u_j$ are adjacent, $|u_j-v_{j-1}|=1$. Hence one of
the first $j-1$ entries of $u_j$ is equal to $x$ and the rest are 1
for some $x\neq 1$, while the last $d-(j-1)$ entries are 0.  Thus
$|u_j-v_0|=j-2$ or $j-1$ according as $x$ is 0 or not. However by
$(H*)$, $|u_j-v_0|=j$ which is a contradiction. Thus, $u_j$ has at
least one non-zero entry in the last $d-(j-1)$ entries. Further,
since the last $d-(j-1)$ entries of $v_{j-1}$ are 0 and
$|u_j-v_{j-1}|=1$, it follows that for some $x\in Z_n\setminus
\{0\}$, the first $(j-1)$ entries of $u_j$  are 1, and $x$ is the
unique non-zero entry in the last $d-(j-1)$ entries. Moreover, since
for each  $i\in [0,j-1]$, the last $d-(j-1)$ entries of $v_i$ are 0,
it follows that $Sym(Z_n\setminus \{0\})\wr Sym([j,d])\leq
A_{v_0,\cdots,v_{j-1}}$. Thus  there exists $\b \in
A_{v_0,\cdots,v_{j-1}}$ such that
$u_j^\b=(1,\cdots,1,1,0,\cdots,0)=v_j$. Therefore, $\mathcal{U}^{\a
\b}=\mathcal{V}$, and hence  $\Gamma$ is $j$-geodesic transitive. By
induction, $\Gamma$ is geodesic transitive.

If $n=2$, then for each vertex $u$, and for any two vertices $v,w$
of $\Gamma(u)$, $|v-w|=2$, that is $v,w$ are not adjacent, so the
girth of $\Gamma$ is not 3. Further, $u_1=(0,0,\cdots,0)$,
$u_2=(1,0,\cdots,0)$, $u_3=(0,1,0,\cdots,0)$ and
$u_4=(1,1,0\cdots,0)$ are four vertices of $\Gamma$ such that
$(u_1,u_2,u_4,u_3,u_1)$ is a 4-cycle, so the girth is 4. Now 2-arcs
and 2-geodesics are the same,  and since  $\diam(\Gamma)=d\geq 3$,
it follows that $\Gamma$ is $2$-arc transitive but not $3$-arc
transitive. If $n\geq 3$, then $(u_1,u_2,u_3,u_1)$ is a triangle
where $u_1=(0,0,\cdots,0)$, $u_2=(1,0,0,\cdots,0)$ and
$u_3=(2,0,0,\cdots,0)$, so $\Gamma$ has girth $3$, that is, $\Gamma$
is arc-transitive but not $2$-arc transitive. $\Box$

\subsection{Odd graphs}

Let $\Omega=[1,2k+1]$ where $k\geq 1$. Then the \emph{Odd graph}
$O_{k+1}$ is the graph whose vertex set is the set of all
$k$-subsets of $\Omega$, and two vertices  are adjacent if and only
if they are disjoint.  By \cite[Section 9.1]{BCN},  $O_{k+1}$ is
distance transitive, its diameter is $k$, valency is $k+1$, and
$\Aut O_{k+1} \cong S_{2k+1}$. By \cite{Biggs-2}, for two vertices
$u,v$ of $O_{k+1}$,

\medskip

$u\in \Gamma_i(v)$ if and only if $|u\cap v|=j$ when $i=2j+1$;
$|u\cap v|=k-j$ when $i=2j$. \ \ \ \ \ \ $(O*)$

\medskip

Note that if $k=1$, then $O_2\cong C_3$ is $s$-arc transitive for
all $s\geq 1$. So we will assume that $k\geq 2$.

\begin{prop}\label{odd-d-not-g}
Let $\Gamma=O_{k+1}$ with $k\geq 2$. Then $\Gamma$  is geodesic
transitive, and is $3$-arc transitive but not $4$-arc transitive.
\end{prop}
{\bf Proof.} Since $\Gamma$ is distance transitive, it follows that
$\Gamma$ is 1-geodesic transitive.

Now, suppose that $\Gamma$ is $(j-1)$-geodesic transitive where
$j\in [2,k]$. Let $\mathcal{V}=(v_0,v_1,\cdots,v_j)$ where
$v_{2i}=[1,k-i]\cup [2k-i+2,2k+1]$ and $v_{2i+1}=[k-i+1,2k-i]$ for
$i\geq 0$.  Then by $(O*)$, $\mathcal{V}$ is a $j$-geodesic. Suppose
that $\mathcal{U}$ is any other $j$-geodesic of $\Gamma$. Then since
$\Gamma$ is $(j-1)$-geodesic transitive, there exists $\a \in
A:=\Aut \Gamma$ such that
$\mathcal{U}^\a=(v_0,v_1,\cdots,v_{j-1},u_j)$ for some $u_j$.

First, suppose that $j=2l$ is even. Let $\Delta_1=[1,k-l+1]$,
$\Delta_2=[k-l+2,k]$, $\Delta_3=[k+1,2k-l+1]$ and
$\Delta_4=[2k-l+2,2k+1]$. Then
$\Omega=\Delta_1\cup\Delta_2\cup\Delta_3\cup\Delta_4$,
$v_{j-1}=\Delta_2 \cup \Delta_3$ and $v_j=(\Delta_1\cup
\Delta_4)\setminus \{k-l+1\}$. Since $v_{j-1}$ and $u_j$ are
adjacent, it follows that $v_{j-1}\cap u_j=\emptyset$, and hence
$u_j \subseteq \Delta_1 \cup \Delta_4$. Since $u_j \in
\Gamma_j(v_0)$ and $j$ is even, by $(O*)$, $|v_0\cap u_j|=k-l$.
Hence $|u_j \setminus v_0|=l$. Since $\Delta_1 \subseteq v_0$, it
follows that $u_j \setminus v_0\subseteq \Delta_4$. Since  $|u_j
\setminus v_0|=l=|\Delta_4|$, it follows that $u_j \setminus v_0=
\Delta_4$, so $\Delta_4  \subseteq u_j$. Hence $u_j=(\Delta_1\cup
\Delta_4)\setminus \{x\}$ for some  $x\in \Delta_1$. Since for each
$m\in [0,j-1]$,  $\Delta_1\subseteq v_m$ when $m$ is even, and
$\Delta_1\cap v_m=\emptyset$ when $m$ is odd, it follows that
$Sym(\Delta_1)\leq A_{v_0,v_1,\cdots,v_{j-1}}$. There exists $\b \in
Sym(\Delta_1)$ such that  $x^\b=k-l+1$. Hence $u_j^\b=v_j$ and so
$\mathcal{U}^{\a \b}=\mathcal{V}$.

Second, suppose that $j=2l+1$ is odd. We let $\Delta_1=[1,k-l]$,
$\Delta_2=[k-l+1,k]$, $\Delta_3=[k+1,2k-l+1]$ and
$\Delta_4=[2k-l+2,2k+1]$. Then
$\Omega=\Delta_1\cup\Delta_2\cup\Delta_3\cup\Delta_4$,
$v_{j-1}=\Delta_1 \cup \Delta_4$ and $v_j=(\Delta_2\cup
\Delta_3)\setminus \{2k-l+1\}$, so $u_j \subseteq \Delta_2\cup
\Delta_3$. Since $v_0\cap \Delta_3=\emptyset$, it follows that
$v_0\cap u_j \subseteq \Delta_2$. Further since $j$ is odd, by
$(O*)$,  $|v_0\cap u_j|=l=|\Delta_2|$. Hence $v_0\cap u_j =
\Delta_2$. Thus $u_{j}=(\Delta_2 \cup \Delta_3)\setminus \{x\}$ for
some  $x\in \Delta_3$. Since for each $m\in [0,j-1]$,
$\Delta_3\subseteq v_m$ when $m$ is odd, and $\Delta_3\cap
v_m=\emptyset$ when $m$ is even, it follows that $Sym(\Delta_3)\leq
A_{v_0,v_1,\cdots,v_{j-1}}$. There exists $\b \in Sym(\Delta_3)$
such that $x^\b=2k-l+1$, and hence   $u_j^\b=v_j$,  and
$\mathcal{U}^{\a \b}=\mathcal{V}$. Thus $\Gamma$ is $j$-geodesic
transitive. Therefore, by induction $\Gamma$ is geodesic transitive.

If  $k=2$, then $\Gamma$ is the Petersen graph which   has girth $5$
and is $3$-arc transitive but not $4$-arc transitive. If $k\geq 3$,
then by \cite[Section 9.1]{BCN},  $\Gamma$ has girth $6$, that is,
3-arcs and 3-geodesics are the same. This together with geodesic
transitivity show that $\Gamma$ is 3-arc transitive. Let  $k=3$ and
$v_0=\{1,2,3\}$, $v_1=\{4,5,6\}$, $v_2=\{1,2,7\}$, $v_3=\{3,4,5\}$,
$v_4=\{1,2,6\}$ and $v_5=\{1,6,7\}$.  Then
$\mathcal{W}_1=(v_0,v_1,v_2,v_3,v_4)$ and
$\mathcal{W}_2=(v_0,v_1,v_2,v_3,v_5)$ are two 4-arcs,
$d_{\Gamma}(v_0,v_4)=2$ and $d_{\Gamma}(v_0,v_5)=3$. So there is no
automorphism mapping $\mathcal{W}_1$ to $\mathcal{W}_2$, and hence
$\Gamma$ is not 4-arc transitive. If $k\geq 4$, then
$\diam(\Gamma)=k\geq 4$ and some 4-arcs lie in 6-cycles and so are
not 4-geodesics. Hence $\Gamma$ is $3$-arc transitive but not
$4$-arc transitive. $\Box$

\subsection{Proof of Theorem \ref{gtnotat}}

First we collect information about the geodesic transitivity of
several 4-arc transitive graphs.

\begin{lemma}\label{biggs-smith-foster}
The Biggs-Smith graph and the Foster graph  have valency and
diameter as in Table 1, and  are geodesic transitive. Moreover, for
$s$ as in Table {\rm 1}, these graphs are  $s$-arc transitive but
not $(s+1)$-arc transitive.
\end{lemma}
{\bf Proof.}  Let  $(\Gamma,s)\in \{(Biggs{-}Smith \ graph,4),$
$(Foster \ graph,5)\}$. Then by \cite[p.221]{BCN} and \cite[Theorem
1.1]{Weiss-2}, $\Gamma$ has valency and diameter as in Table 1, and
for $s$ as in Table 1, $\Gamma$ is $s$-arc transitive but not
$(s+1)$-arc transitive. Thus $\Gamma$ is $s$-geodesic transitive.
Let $d=\diam(\Gamma)$ and $(v_0,v_1,\cdots,v_d)$ be a $d$-geodesic.
Then $d=s+3$. By \cite[p.221]{BCN}, $|\Gamma(v_j)\cap
\Gamma_{j+1}(v_0)|= 1$ for every $j=d-3,d-2,d-1$,  and it follows
that $\Gamma$ is geodesic transitive. $\Box$

\bigskip

As in \cite[p.84]{GR} we define a \emph{generalized polygon}, or
more precisely, a \emph{generalized $d$-gon}, as a bipartite graph
with diameter $d$ and girth $2d$. The generalized polygons related
to the Lie type groups $A_2(q),B_2(q)$ and $G_2(q)$ ($q$ is prime
power) are \emph{classical generalized polygons}, and  are denoted
by $\Delta_{3,q}$, $\Delta_{4,q}$ and $\Delta_{6,q}$, respectively.
They  are regular of valency $q+1$.

\begin{lemma}\label{strans-gp-gt}
The only distance transitive generalized polygons of valency at
least $3$ that are $s$-arc transitive but not $(s+1)$-arc
transitive, for some $s\geq 4$, are $\Delta_{s-1,q}$ where $(s,q)\in
S=\{(4,q),(5,2^m),(7,3^m)|$ $q$ is a prime power and $m$ is a
positive integer $\}$. Moreover,  all these graphs are geodesic
transitive.
\end{lemma}
{\bf Proof.} Let $\Gamma$ be a distance transitive  generalized
polygon of valency at least 3 that is $s$-arc transitive but not
$(s+1)$-arc transitive for some $s\geq 4$. Let  $g$ be its girth.
Suppose that $s<\frac{g-2}{2}$. Then for any two vertices $u,v$ at
distance  $s+1$, there exists a unique $(s+1)$-arc between them.
Since $\Gamma$ is distance transitive, it follows that $\Gamma$ is
$(s+1)$-arc transitive, which contradicts our assumption. Thus $s
\geq \frac{g-2}{2}$. Since $\Gamma$  is a generalized polygon, $g$
is even. By \cite[Theorem 1.1]{Weiss-2}, $s\in \{ \frac{g+2}{2},
\frac{g}{2},\frac{g-2}{2} \}$.   If $s= \frac{g+2}{2}$ or
$\frac{g}{2}$, then  \cite[Theorem 1.1]{Weiss-2} shows that $\Gamma$
is one of $\Delta_{s-1,q}$ where $(s,q)\in S$. Let $A=\Aut \Gamma$.
Since $\Gamma$ is distance transitive, $A_u$ is transitive on
$\Gamma_{s+1}(u)$ for each vertex $u$.  Thus, if $s= \frac{g-2}{2}$,
by \cite[Theorem 1.1]{Weiss-2}, $\Gamma$ is the Biggs-Smith graph,
which is not a generalized polygon.  Moreover, since each
$\Delta_{s-1,q}$ has  diameter  $s-1$, it follows that  all these
graphs are geodesic transitive. $\Box$

\medskip
These lemmas allow us to specify precisely the geodesic transitive
graphs which are 4-arc transitive.

\begin{prop}\label{dt-4at-diam}
Let $\Gamma$ be a regular graph of valency at least $3$. Then
$\Gamma$ is geodesic transitive and  $s$-arc transitive but not
$(s+1)$-arc transitive for some $s\geq 4$ if and only if $\Gamma$
is in one of the lines of Table 1. Table 1 also gives the valency,
 integer $s$ and diameter for each graph.   In particular, $\diam(\Gamma)\leq 8$.
\begin{table}[!hbp]\caption{Geodesic and  $s$-arc transitive graphs that  not $(s+1)$-arc transitive  }
\medskip
\centering
\begin{tabular}{|c|c|c|c|c|}
\hline
 Graph $ \Gamma$  &  Valency & $s$ & Diameter \\
\hline
Foster graph  & $3$ & 5& 8 \\
\hline
Biggs-Smith graph   & $3$ & 4 & 7  \\
\hline
   $\Delta_{3,q}$, $q$ is a prime power   & $q+1$ & 4& 3 \\
 \hline
$\Delta_{4,q}$, $q=2^m$, $m$ is a positive integer    &  $q+1$ & 5&4\\
\hline
$\Delta_{6,q}$, $q=3^m$, $m$ is a positive integer    & $q+1$ &  7&6 \\
  \hline

\end{tabular}

\end{table}

\end{prop}
{\bf Proof.}  By Corollary 1.2 of \cite{Weiss-2},  if $\Gamma$  is
distance transitive, $s$-arc transitive but not $(s+1)$-arc
transitive for some $s\geq 4$, then $\Gamma$ is either the Foster
graph or the Biggs-Smith graph, or a generalized polygon. The result
now follows from  Lemmas \ref{biggs-smith-foster} and
\ref{strans-gp-gt}.  $\Box$

\bigskip

\noindent {\bf Proof of Theorem \ref{gtnotat}.} It follows from
Propositions \ref{j(n,k)}, \ref{Ham-2-trans}, \ref{odd-d-not-g} and
\ref{dt-4at-diam} that, for each  $s\in \{1,2,3,4,5,7\}$, there are
infinitely many geodesic transitive graphs that are $s$-arc
transitive but not $(s+1)$-arc transitive. Further, for each $s\in
\{1,2,3\}$, Propositions \ref{j(n,k)}, \ref{Ham-2-trans} and
\ref{odd-d-not-g} show that  there are such  graphs  with
arbitrarily large diameter. Finally, by Proposition
\ref{dt-4at-diam}, this is not the case for  $s\in \{4,5,7\}$.
Therefore geodesic transitive graphs that are $s$-arc transitive but
not $(s+1)$-arc transitive with arbitrarily large diameter occur
only for $s\in \{1,2,3\}$. $\Box$

\section{ Paley graphs}

In this section, we  discuss a special family of connected Cayley
graphs, namely the Paley graphs, which  were first defined by Paley
in 1933, see \cite{Paley-1}. We   prove that the Paley graph $P(q)$
is distance transitive but not geodesic transitive whenever $q\geq
13$.

Let $q=p^e$ be a prime power such that $q\equiv 1 \pmod{4}$. Let
$F_q$ be a finite field of order $q$. The  \emph{Paley graph} $P(q)$
is the graph with vertex set  $F_q$, and two distinct vertices $u,v$
are adjacent if and only if $u-v$ is a nonzero square  in $F_q$. The
congruence condition on $q$ implies that $-1$ is a square in $F_q$,
and hence $P(q)$ is an undirected graph.

Note that the field $F_q$ has $\frac{q-1}{2}$ elements which are
nonzero squares, so $P(q)$ has valency $\frac{q-1}{2}$. Moreover,
$P(q)$ is a Cayley graph for the additive group $G=F_{q}^+\cong
Z_{p}^e$. Let $w$ be a primitive element of $F_q$.  Then
$S=\{w^2,w^4,\cdots,w^{q-1}=1\}$ is the set of  nonzero squares  of
$F_q$,  and $P(q)=\Cay(G,S)$. Define $\tau:F_q\mapsto F_q,x\mapsto
x^p$. Then $\tau$ is an automorphism of the field $F_q$, called the
\emph{Frobenius automorphism}, and $\Aut F_q=\langle \tau
\rangle\cong Z_e$.  By \cite{Carlitz} (or see \cite{Lim-Praeger-1}),
$\Aut P(q)=(G:\langle w^2\rangle).\langle \tau \rangle\cong
(Z_p^e:Z_{\frac{q-1}{2}}).Z_e\leq \A\Gamma L(1,q)$.

Let  $S'$ be the set of all  nonsquare elements of $G$. Then
$|S'|=\frac{q-1}{2}$.  Define the Cayley graph $\Sigma=\Cay(G,S')$
where two vertices $u,v$ are adjacent if and only if $u-v\in S'$.
Then $\Sigma$ is the complement of the Paley graph $P(q)$. Further,
multiplication by $w$ induces an isomorphism $\Sigma\cong P(q)$, see
\cite{Sachs}.

Now, we cite a property of Paley graphs.

\begin{lemma}{\rm(\cite[p.221]{GR})}\label{val-p-rem1}
Let  $\Gamma=P(q)$, where  $q$ is a prime power such that $q\equiv 1
\pmod 4$. Let $u,v$ be   distinct vertices of $\Gamma$. If $u,v$ are
adjacent, then $|\Gamma(u)\cap \Gamma(v)|=\frac{q-5}{4}$; if $u,v$
are not adjacent, then $|\Gamma(u)\cap \Gamma(v)|=\frac{q-1}{4}$.
\end{lemma}

\subsection{Proof of Theorem \ref{exam-valp-ha}}

Let $F_q, G$ and $w$ be as above, and let $G^*=\langle w \rangle $
be the multiplicative group  of $F_q$.  As discussed above,
$P(q)=\Cay(G,S)$ where $S=\{w^2,w^4,\cdots,w^{q-1}=1\}$, and  $\Aut
P(q)=(G:\langle w^2\rangle).\langle \tau \rangle\cong
(Z_p^e:Z_{\frac{q-1}{2}}).Z_e$.

Now we  prove Theorem \ref{exam-valp-ha}.

\bigskip

\noindent {\bf Proof of Theorem \ref{exam-valp-ha}.} Let
$\Gamma=P(q)$ and $A=\Aut \Gamma$.  Let $u=0 \in G$. Then
$A_{u}=\langle w^2 \rangle.\langle \tau \rangle$ has orbits $\{0\}$,
$S$ and $S'=G\setminus (\{0\}\cup S)$ on vertices. Now $S=\Gamma(u)$
and as $\Gamma$ is connected, $\Gamma_2(u)$ must be the other orbit
$S'$. In particular, $\Gamma$ has diameter 2 and is distance
transitive and arc transitive.

Suppose that $v\in \Gamma(u)$. By Lemma \ref{val-p-rem1},
$|\Gamma(u)\cap \Gamma(v)|=\frac{q-5}{4}$. Thus, $|\Gamma_2(u)\cap
\Gamma(v)|=|\Gamma(v)|-|\Gamma(u)\cap \Gamma(v)|-1=\frac{q-1}{4}$.
If $A$ is transitive on the 2-geodesics of $\Gamma$ then $A_{u,v}$
is transitive on $\Gamma_2(u)\cap \Gamma(v)$. In particular
$\frac{q-1}{4}=|\Gamma_2(u)\cap \Gamma(v)|$ divides $|A_{uv}|=e$ and
hence $q=5$ or 9. We consider $P(5)$ and $P(9)$.

If $q=5$, then  $P(q)\cong C_5$, so $P(q)$ is geodesic transitive.

Now, suppose that $q=9=3^2$. The field $F_9$ is $\{a+bx \,|\, a,b
\in Z_3\}$ under polynomial addition and multiplication modulo
$f(x)=x^2+1$.  The set $S$ is $\{1,2,x,2x\}$, and
$\Gamma_2(u)=\{x+1,x+2,2x+1,2x+2 \}$. Let $v=1\in S$. Then
$\Gamma_2(u)\cap \Gamma(v)=\{x+1,2x+1\}$. Since  $A_{uv}=\langle
\tau \rangle$  and $(x+1)^{\tau}=x^3+1=2x+1$, it follows that
$A_{uv}$ is transitive on $\Gamma_2(u)\cap \Gamma(v)$, and hence
$\Gamma$ is $(A,2)$-geodesic transitive. Since $\diam(\Gamma)=2$,
$\Gamma=P(9)$ is geodesic transitive. $\Box$

\subsection{Arc-transitive graphs of odd prime order}

In this subsection we characterise the Paley graphs $P(p)$, for
primes $p$, as arc-transitive graphs of given prime order and given
valency. This result is used in our proof of Theorem
\ref{gt-primeval}.

\begin{prop}\label{arctrans-order-p-1}
Let $\Gamma$ be an arc-transitive graph of  prime order $p$ and
valency $\frac{p-1}{2}$. Then $p\equiv 1 \pmod 4$, $\Aut \Gamma
\cong Z_p:Z_{\frac{p-1}{2}}$, and  $\Gamma \cong P(p)$.
\end{prop}

The proof uses the following famous result of Burnside.

\begin{lemma}{\rm( \cite[Theorem 3.5B]{DM-1})}\label{val-p-prim-1}
Suppose that $G$ is a primitive permutation group  of prime degree
$p$. Then $G$ is either $2$-transitive, or solvable and $G\leq
AGL(1,p)$.
\end{lemma}

\noindent {\bf Proof of Proposition \ref{arctrans-order-p-1}.} Since
$\Gamma$ has  valency $\frac{p-1}{2}$,  $p$ is an odd prime. Since
$\Gamma$ is  undirected and arc-transitive, it follows that $\Gamma$
has $p(\frac{p-1}{2})/2$ edges. This implies that  $p\equiv 1 \pmod
4$.

Let $A=\Aut \Gamma$. Since  $A$ is transitive on $V\Gamma$ and  $p$
is a prime, $A$ is primitive on $V\Gamma$. Since $\Gamma$ is neither
complete nor empty, it follows by Lemma \ref{val-p-prim-1} that $A<
AGL(1,p)=Z_p:Z_{p-1}$. Again by vertex transitivity,  $Z_p\leq A$.
Thus, $A \cong Z_p:Z_m$ where $Z_m< Z_{p-1}$.

Since $Z_p$ is regular on $V\Gamma$, it follows from Lemma
\ref{cayley-1} that $\Gamma$ is a Cayley graph for $Z_p$. Thus
$\Gamma=\Cay(G,S)$ where $G\cong Z_p$, $S\subseteq G\setminus
\{0\}$, $S=S^{-1}$ and $|S|=\frac{p-1}{2}$. Let $v\in V\Gamma$ be
the vertex corresponding to $0\in G$. Then $A_v=Z_m$ acts
semiregularly on $G\setminus \{v\}$ with orbits of size $m$. Since
$\Gamma$ is arc-transitive, $A_v$ acts transitively on $S$, so
$m=|S|=\frac{p-1}{2}$. Thus $A \cong Z_p:Z_{\frac{p-1}{2}}$.

Now we may identify $G$ with $F_p^+$ and $v$ with 0. Then $A_v$ is
the unique subgroup of order $\frac{p-1}{2}$ of $F_p^*=\langle
w\rangle$, that is, $A_v=\langle w^2\rangle$. The $A_v$-orbits in
$F_p$ are $\{0\}$, $S_1=\{w^2,w^4,\cdots,w^{p-1}\}$ and
$S_2=\{w,w^3,\cdots,w^{p-2}\}$, and so $S=S_1$ or $S_2$, and
$\Gamma=P(p)$ or its complement respectively. In either case,
$\Gamma\cong P(p)$.  $\Box$

\section{ Graphs of prime valency that are 2-geodesic transitive but not 2-arc transitive}

We  prove Theorem \ref{gt-primeval} in Subsection 5.1, that is, we
give a classification of connected 2-geodesic transitive graphs of
prime valency which are not 2-arc transitive. Note that the
assumption of 2-geodesic transitivity implies that the graph is not
complete. Since the identification of the examples is made by
reference to deep classification results, we give in Subsection 5.2
an explicit  construction of these graphs as coset graphs and verify
the properties claimed in Theorem \ref{gt-primeval}.

\subsection{ Proof of Theorem \ref{gt-primeval}}

We will prove Theorem \ref{gt-primeval} in a series of lemmas.
Throughout this subsection we assume that $\Gamma$ is a connected
$2$-geodesic transitive but not $2$-arc transitive  graph of prime
valency $p$ and that $A=\Aut \Gamma$.  The first Lemma
\ref{val-p-lemma-1} determines some intersection parameters.

\begin{lemma}\label{val-p-lemma-1}
Let $(v,u,w)$ be a $2$-geodesic  of $\Gamma$. Then $p\equiv 1 \pmod
4$, $|\Gamma(v)\cap \Gamma(u)|=|\Gamma_2(v)\cap
\Gamma(u)|=\frac{p-1}{2}$ and $|\Gamma(v)\cap \Gamma(w)|$ divides
$\frac{p-1}{2}$. Moreover, $A_v^{\Gamma(v)} \cong
Z_p:Z_{\frac{p-1}{2}}$, $A_{v,u}^{\Gamma(v)} \cong
Z_{\frac{p-1}{2}}$ and $A_{v,u}$ is transitive on $\Gamma(v)\cap
\Gamma(u)$.
\end{lemma}
{\bf Proof.} Since $\Gamma$ is $2$-geodesic transitive but not 2-arc
transitive, it follows that  $\Gamma$ is not a cycle. In particular,
$p$ is an odd prime. Let  $|\Gamma(v)\cap \Gamma(u)|=x$ and
$|\Gamma_2(v)\cap \Gamma(u)|=y$. Then $x+y=|\Gamma(u)\setminus
\{v\}|=p-1$. Since $\Gamma$ is 2-geodesic transitive but not 2-arc
transitive, it follows that  $\girth(\Gamma)=3$, so $x\geq 1$. Since
the induced subgraph $[\Gamma(v)]$ is an undirected regular graph
with $\frac{px}{2}$ edges, and since $p$ is  odd, it follows that
$x$ is even. This together with  $x+y=p-1$ and the fact that $p-1$
is even, implies that $y$ is also even.

Since $\Gamma$ is arc-transitive, $A_v^{\Gamma(v)}$ is transitive on
$\Gamma(v)$.  Since $p$ is a prime, $A_v^{\Gamma(v)}$ acts
primitively on $\Gamma(v)$. By Lemma \ref{val-p-prim-1}, either
$A_v^{\Gamma(v)}$ is 2-transitive, or $A_v^{\Gamma(v)}$ is solvable
and $A_v^{\Gamma(v)}\leq AGL(1,p)$. Since $\Gamma$ is not complete,
it follows that $[\Gamma(v)]$ is not a  complete graph. Also since
$\girth(\Gamma)=3$, $[\Gamma(v)]$ is not an empty graph and so
$A_v^{\Gamma(v)}$ is not 2-transitive. Hence $A_v^{\Gamma(v)}<
AGL(1,p)$. Thus $A_v^{\Gamma(v)}\cong Z_p:Z_m$, where $m|(p-1)$ and
$m<p-1$. Hence $m\leq \frac{p-1}{2}$.

Since $\Gamma$ is vertex transitive, it follows that
$A_u^{\Gamma(u)}\cong Z_p:Z_m$, and hence $A_{u,v}^{\Gamma(u)}\cong
Z_m$ is semiregular on $\Gamma(u)\setminus \{v\}$ with orbits of
size $m$. Since  $\Gamma$ is $2$-geodesic transitive,
$A_{u,v}^{\Gamma(u)}$ is transitive on $\Gamma_2(v)\cap \Gamma(u)$,
and hence  $y=|\Gamma_2(v)\cap \Gamma(u)|=m$, so
$x=p-1-m=m(\frac{p-1}{m}-1)\geq m$, and $x$ is divisible by $m$.

Now  again by arc transitivity, $|\Gamma(u)\cap
\Gamma(w)|=|\Gamma(u)\cap \Gamma(v)|=x$. Since $|\Gamma_2(v)\cap
\Gamma(u)|=m$, it follows that $|\Gamma_2(v)\cap \Gamma(u)\cap
\Gamma(w)|\leq m-1$. Since $\Gamma(w)\cap \Gamma(u)=(\Gamma(w)\cap
\Gamma(u)\cap \Gamma(v))\cup (\Gamma(w)\cap \Gamma(u)\cap
\Gamma_2(v))$,  it follows that
$$x\leq |\Gamma(w)\cap \Gamma(u)\cap \Gamma(v)|+(m-1). \ \ \ \ \ \ \  (*)$$

Let $z=|\Gamma(v)\cap \Gamma(w)|$ and $n=|\Gamma_2(v)|$. Since
$\Gamma$ is $2$-geodesic transitive, $z,n$ are independent of $v,w$
and, counting edges between $\Gamma(v)$ and $\Gamma_2(v)$ we have
$pm=nz$. Now $z\leq |\Gamma(v)|=p$. Suppose first that $z=p$. Then
$m=n$ and $\Gamma(v)=\Gamma(w)$, and so for distinct $w_1,w_2\in
\Gamma_2(v)$, $d_{\Gamma}(w_1,w_2)=2$. Since $\Gamma$ is 2-geodesic
transitive, it follows that $\Gamma(v)=\Gamma(v')$ whenever
$d_{\Gamma}(v,v')=2$. Thus $\diam(\Gamma)=2$, $V\Gamma=\{v\}\cup
\Gamma(v)\cup \Gamma_2(v)$ and $|V\Gamma|=1+p+m$. Let
$\Delta=\{v\}\cup \Gamma_2(v)$. Then for distinct $v_1,v_1'\in
\Delta $, $d_{\Gamma}(v_1,v_1')=2$; for any $v_1''\in
V\Gamma\setminus \Delta$, $v_1,v_1''$ are adjacent. Thus, for any
$v_1\in \Delta$, $\Delta=\{v_1\}\cup \Gamma_2(v_1)$. It follows that
$\Delta$ is a block of imprimitivity for  $A$ of size $m+1$. Hence
$(m+1)|(p+m+1)$, so $(m+1)|p$. Since $m|(p-1)$,  it follows that
$m+1=p$ which contradicts the inequality $m\leq \frac{p-1}{2}$.

Thus $z<p$, and so  $z$ divides $m$.  Since $|\Gamma(w)\cap
\Gamma(u)\cap \Gamma(v)|\leq z$, it follows from  $(*)$ that $x\leq
z+(m-1)\leq 2m-1<2m$.  Since $x$ is divisible by $m$ and $x\geq m$
we have $x=m$. Thus $2m=x+y=p-1$, so $x=y=m=\frac{p-1}{2}$, and
since $x$ is even, $p\equiv 1 \pmod 4$. Also $x=m$ implies that
$A_{v,u}$ is transitive on $\Gamma(v)\cap \Gamma(u)$. Finally, since
$nz=pm=p(\frac{p-1}{2})$ and $z<p$, it follows that $z$ divides
$\frac{p-1}{2}$.  $\Box$

\begin{lemma}\label{val-p-lemma-2}
For $v\in V\Gamma$,  the stabiliser $A_v \cong
Z_p:Z_{\frac{p-1}{2}}$.
\end{lemma}
{\bf Proof.} Suppose that  $(v,u)$ is an arc of $\Gamma$.  Then by
Lemma \ref{val-p-lemma-1}, $A_v^{\Gamma(v)}\cong
Z_p:Z_{\frac{p-1}{2}}$, and $A_{v,u}^{\Gamma(v)} \cong
Z_{\frac{p-1}{2}}$ is regular on $\Gamma(v)\cap \Gamma(u)$. Let $E$
be  the kernel of the action of $A_v$ on $\Gamma(v)$.  Let $u'\in
\Gamma(v)\cap \Gamma(u)$ and $x\in E$. Then $x\in A_{v,u,u'}$. Since
$A_{u,v}^{\Gamma(u)} \cong Z_{\frac{p-1}{2}}$ is semiregular on
$\Gamma(u)\setminus \{v\}$, it follows that $x$ fixes all vertices
of $\Gamma(u)$. Since $x$ also fixes all vertices of $\Gamma(v)$,
this argument for each $u\in \Gamma(v)$ shows that  $x$ fixes all
vertices of $\Gamma_2(v)$. Since $\Gamma$ is connected, $x$ fixes
all vertices of $\Gamma$, hence $x=1$. Thus $E=1$, so $A_v \cong
Z_p:Z_{\frac{p-1}{2}}$. $\Box$

\begin{lemma}\label{val-p-lemma-3}
Let $(v,u,w)$ be a $2$-geodesic of $\Gamma$. Then $|\Gamma(v)\cap
\Gamma(w)|=\frac{p-1}{2}$, $|\Gamma_2(v)\cap \Gamma(w)\cap
\Gamma(u)|=\frac{p-1}{4}$,  $|\Gamma_2(v)|=p$, and $|\Gamma_2(v)\cap
\Gamma(w)|=\frac{p-1}{2}$.
\end{lemma}
{\bf Proof.}  Let $z=|\Gamma(v)\cap \Gamma(w)|$ and
$n=|\Gamma_2(v)|$. By  Lemma \ref{val-p-lemma-1}, $|\Gamma(u)\cap
\Gamma_2(v)|=\frac{p-1}{2}$ and $z|\frac{p-1}{2}$. Counting the
edges between $\Gamma(v)$ and $\Gamma_2(v)$ gives
$\frac{p-1}{2}p=nz$. By Lemma \ref{val-p-lemma-2},
$A_{v,u}=Z_{\frac{p-1}{2}}$, and by Lemma \ref{val-p-lemma-1},
$A_{v,u}$ is transitive on $\Gamma(v)\cap \Gamma(u)$, so
$[\Gamma(u)]$ is $A_{u}$-arc transitive. Since $p$ is a prime, it
follows by Lemma \ref{arctrans-order-p-1} that $[\Gamma(u)]$ is a
Paley graph $P(p)$. Since $v,w\in \Gamma(u)$ are not adjacent, by
Lemma \ref{val-p-rem1}, $|\Gamma(v)\cap \Gamma(u)\cap
\Gamma(w)|=\frac{p-1}{4}$, hence $z\geq \frac{p-1}{4}+1$. Since
$z|\frac{p-1}{2}$, it follows that $z=\frac{p-1}{2}$. Hence $n=p$.
Thus, $|\Gamma(v)\cap \Gamma(w)|=\frac{p-1}{2}$ and
$|\Gamma_2(v)|=p$.

By Lemma \ref{val-p-lemma-1}, we have $|\Gamma(v)\cap
\Gamma(u)|=\frac{p-1}{2}$. Since $\Gamma$ is arc transitive, it
follows that  $|\Gamma(v_1)\cap \Gamma(v_2)|=\frac{p-1}{2}$ for
every arc $(v_1,v_2)$. Thus, $|\Gamma(u)\cap
\Gamma(w)|=\frac{p-1}{2}$. Since $\Gamma(u)\cap
\Gamma(w)=(\Gamma(v)\cap \Gamma(u)\cap \Gamma(w) )\cup
(\Gamma_2(v)\cap \Gamma(u)\cap \Gamma(w) )$ and $|\Gamma(v)\cap
\Gamma(u)\cap \Gamma(w)|=\frac{p-1}{4}$, it follows that
$|\Gamma_2(v)\cap \Gamma(u)\cap
\Gamma(w)|=\frac{p-1}{2}-\frac{p-1}{4}=\frac{p-1}{4}$.  Since $A_v=
Z_p:Z_{\frac{p-1}{2}}$, it follows that $A_{v,w}=Z_{\frac{p-1}{2}}$
and $A_{v,w}$  is semiregular on $\Gamma_2(v)\setminus \{w\}$ with
orbits of size $\frac{p-1}{2}$. Since $\Gamma_2(v)\cap
\Gamma(w)\subseteq \Gamma(w)\setminus \Gamma(v)$ ( of size
$\frac{p-1}{2}$ ) and since  $|\Gamma_2(v)\cap \Gamma(w)\cap
\Gamma(u)|=\frac{p-1}{4}>0$,  it follows that $|\Gamma_2(v)\cap
\Gamma(w)|=\frac{p-1}{2}$. $\Box$

\begin{lemma}\label{val-p-lemma-4}
Let $v$ be a vertex of $\Gamma$. Then $|\Gamma_3(v)|=1$ and
$\diam(\Gamma)=3$. Further, $\Gamma$ is geodesic transitive.
\end{lemma}
{\bf Proof.} Suppose that $(v,u,w)$ is a $2$-geodesic of $\Gamma$.
Then by Lemma \ref{val-p-lemma-3}, $|\Gamma(v)\cap
\Gamma(w)|=\frac{p-1}{2}$ and $|\Gamma_2(v)\cap
\Gamma(w)|=\frac{p-1}{2}$. Hence $|\Gamma_3(v)\cap
\Gamma(w)|=p-|\Gamma(v)\cap \Gamma(w)|-|\Gamma_2(v)\cap
\Gamma(w)|=1$. Since $\Gamma$ is $2$-geodesic transitive, it follows
that $|\Gamma_3(v)\cap \Gamma(w_1)|=1$ for all $w_1\in \Gamma_2(v)$.
Further $\Gamma$ is $3$-geodesic transitive.

Let $\Gamma_3(v)\cap \Gamma(w)=\{v'\}$,  $n=|\Gamma_3(v)|$ and
$i=|\Gamma_2(v)\cap \Gamma(v')|$. Counting edges between
$\Gamma_2(v)$ and  $\Gamma_3(v)$,  we have  $p=ni$.  Since
$[\Gamma(w)]$ is a Paley graph and $u,v'\in \Gamma(w)$ are not
adjacent, it follows from Lemma \ref{val-p-rem1} that
$|\Gamma(u)\cap \Gamma(w)\cap \Gamma(v')|=\frac{p-1}{4}$. Since
$\Gamma(u)\cap \Gamma_2(v)$ contains these $\frac{p-1}{4}$ vertices
as well as $w$, we have $i\geq \frac{p+3}{4}>1$. Thus  $i=p$ and
$n=1$, that is, $|\Gamma_3(v)|=1$. Since $|\Gamma_2(v)\cap
\Gamma(v')|=p$ and $|\Gamma_2(v)|=p$, it follows that
$\Gamma_2(v)=\Gamma(v')$, and so  $\diam(\Gamma)=3$. Therefore
$\Gamma$ is geodesic transitive. $\Box$

\bigskip

Now, we  prove Theorem \ref{gt-primeval}.

\bigskip

\noindent {\bf Proof of Theorem \ref{gt-primeval}.} Since the graph
$\Gamma$ is $2$-geodesic transitive but not $2$-arc transitive, it
follows that $\girth(\Gamma)=3$, and hence $\Gamma$ is nonbipartite.
Let $v\in V\Gamma$. Then it follows from Lemmas \ref{val-p-lemma-1}
to \ref{val-p-lemma-4} that   $p\equiv 1 \pmod 4$,
$|\Gamma_2(v)|=p$, $|\Gamma_3(v)|=1$ and $\diam(\Gamma)=3$. Thus,
$V\Gamma=\{v\}\cup \Gamma(v)\cup \Gamma_2(v)\cup \{v'\}$, where
$\Gamma_3(v)=\{v'\}$, $\Gamma(v)=\Gamma_2(v')$ and
$\Gamma_2(v)=\Gamma(v')$.    Since $\Gamma$ is vertex transitive,
these properties hold  for all vertices of $\Gamma$. Thus, $\Gamma$
is an antipodal graph.  By Lemma \ref{val-p-lemma-4}, $\Gamma$ is
geodesic transitive,  and hence distance transitive.

Let $\mathcal{B}=\{\Delta_1,\Delta_2, \cdots,\Delta_{p+1}\}$ where
$\Delta_i=\{u_i,u_i'\}$ such that $d_{\Gamma}(u_i,u_i')=3$. Then
each $\Delta_i$ is a block for $\Aut \Gamma$ of size 2 on $V\Gamma$.
Further, for each $j\neq i$, $u_i$ is adjacent to exactly one vertex
of $\Delta_j$, and $u_i'$ is adjacent to the other. The quotient
graph  $\Sigma$  such that $V\Sigma=\mathcal{B}$, and two vertices
$\Delta_i,\Delta_j$ are adjacent if and only if
$\{\Delta_i,\Delta_j\}$ contains an edge of $\Gamma$, is therefore a
complete graph  $\Sigma\cong K_{p+1}$ and $\Gamma$ is a cover of
$\Sigma$.

Therefore, we know that  $\Gamma$ is a nonbipartite antipodal
distance transitive cover with fibres of size 2  of the complete
graph $K_{n}$, where $n=p+1$, $p\equiv  1\pmod 4$, and
$\diam(\Gamma)=3$, so it is one of the graphs listed in I or II of
Lemma \ref{anticover-lemma-1}. Since $p\equiv 1\pmod 4$, it follows
that $n\equiv 2\pmod 4$. However, for I $(a)$,  $(c)$, and II,
$n\equiv 0\pmod 4$,  so $\Gamma$ is not a graph in one of these
cases. For I $(b)$ and $(d)$, $n-1$ is not a prime, so $\Gamma$ is
not in one of these cases either. Thus $\Gamma$ is the graph in I
$(e)$ of Lemma \ref{anticover-lemma-1} with $q=p$ prime. Hence
$\Gamma$ is unique up to isomorphism and $A=\Aut(\Gamma)\leq
PSL(2,p) \times Z_2$. By Lemma \ref{val-p-lemma-2},
$A_v=Z_p:Z_{\frac{p-1}{2}}$ for every $v\in V\Gamma$, and since
$\Gamma$ is vertex transitive, it follows that
$|A|=p(p+1)(p-1)=|PSL(2,p)\times Z_2|$. Thus $A=  PSL(2,p)\times
Z_2$. $\Box$

\subsection{Construction }

In Subsection 5.1, we  proved  Theorem \ref{gt-primeval}  using
results of \cite{GLP} to   classify  the connected 2-geodesic
transitive graphs of prime valency $p$ which are not 2-arc
transitive. Here we identify the examples explicitly. The unique
example of valency 5 is the icosahedron, and we assume from now on
that $p>5$  and $p\equiv 1\pmod{4}$. Taylor gave a construction of
this family of graphs from regular two-graphs, see \cite[p.14]{BCN}
and \cite{Taylor-1}. Here  we present a direct construction of these
graphs as coset graphs, fleshing out the  construction given by the
third author in the proof of \cite[Theorem 1.1]{Li-ci-soluble} in
order to prove the additional properties we need for  Theorem
\ref{gt-primeval}.

For a finite group $G$,  a  core-free proper subgroup $H$, and  an
element $g\in G$ such that $G=\langle H,g\rangle$ and $g^2\in H$,
the \emph{coset graph} $\Cos(G,H,HgH)$ is the  graph with vertex set
$\{Hx|x\in G\}$, and two vertices $Hx,Hy$  adjacent if and only if
$yx^{-1}\in HgH$. It is  a  connected, undirected, and  $G$-arc
transitive  graph of valency $|H:H\cap H^g|$, see \cite{Lorimer-1}.

\begin{constr}\label{coset-constr}
{\rm Let $G=PSL(2,p)$ where $p>5$ is a prime  and  $p\equiv
1\pmod{4}$. Choose $a\in G$ such that  $o(a)=p$. Then $N_{G}(\langle
a\rangle)=\langle a\rangle:\langle b\rangle\cong
Z_p:Z_{\frac{p-1}{2}}$ for some  $b\in G$, $o(b)=\frac{p-1}{2}$.
Further, there exists an involution $g\in G$ such that
$N_{G}(\langle b^2\rangle)=\langle b\rangle:\langle g\rangle\cong
D_{p-1}$. Let $H=\langle a\rangle:\langle b^2\rangle$ and
$\Gamma=\Cos(G,H,HgH)$.   }
\end{constr}

First, in the following lemma, we show that the coset graph in
Construction \ref{coset-constr}  is unique up to isomorphism for
each $p$. We repeatedly use the fact that each $\sigma \in \Aut G$
induces an isomorphism from $\Cos(G,H,HgH)$ to
$\Cos(G,H^{\sigma},H^{\sigma}g^{\sigma}H^{\sigma})$, and in
particular, we use this fact for the conjugation action by elements
of $G$.

\begin{lemma}\label{coset-iso-rem}
For each fixed prime $p>5$ and  $p\equiv 1\pmod{4}$, up to
isomorphism, the graph $\Gamma$ in Construction \ref{coset-constr}
is independent of the choices of  $H$ and $g$.
\end{lemma}
{\bf Proof.}  Let $G=PSL(2,p)$ where $p>5$ is a prime  and  $p\equiv
1\pmod{4}$. Let elements  $a_i,b_i,g_i$ and  subgroup $H_i$  be
chosen as  in Construction \ref{coset-constr} for $i\in \{1,2\}$.
Let $X=PGL(2,p)\cong \Aut(G)$.

Since all subgroups of $G$ of order $p$ are conjugate there exists
$x\in G$ such that $\langle a_2\rangle^x=\langle a_1\rangle$, so we
may assume that $\langle a_1\rangle=\langle a_2\rangle=K$, say. Let
$Y=N_{X}(K)$. Then $Y=K: \langle y\rangle$ where $o(y)=p-1$, and
$H_1=K:\langle b_1^2\rangle$ and $H_2=K:\langle b_2^2\rangle$ are
equal to the unique subgroup of $Y$ of order $\frac{p(p-1)}{4}$,
that is, $H_1=H_2=K:\langle y^4\rangle=H$, say. Next, since all
subgroups of $Y$ of order $\frac{p-1}{4}$ are conjugate, there exist
$x_1,x_2\in Y$ such that $\langle b_1^2\rangle^{x_1}=\langle
b_2^2\rangle^{x_2}=\langle y^4\rangle$. Since each $x_i$ normalises
$H$ we may assume in addition that $\langle b_1^2\rangle=\langle
b_2^2\rangle=\langle y^4\rangle < \langle y\rangle$. Thus $g_1,g_2$
are non-central involutions in $N_G(\langle y^4\rangle)\cong
D_{p-1}$, an index 2 subgroup of $N_X(\langle y^4\rangle)=\langle
y\rangle:\langle z\rangle\cong D_{2(p-1)}$. The set of non-central
involutions in $N_G(\langle y^4\rangle)$ form a conjugacy class of
$N_X(\langle y^4\rangle)$ of size $\frac{p-1}{2}$ and consists of
the elements $y^{2i}z$, for $0\leq i<\frac{p-1}{2}$. The group
$\langle y\rangle$ acts transitively on this set of involutions by
conjugation (and normalises $H$). Hence, for some $u\in \langle
y\rangle$, $H^u=H$ and $g_2^u=g_1$.   $\Box$

\bigskip

Now we  show that the coset graph $\Gamma$ in Construction
\ref{coset-constr} is 2-geodesic transitive but not 2-arc transitive
of prime valency $p$. We first state some properties of $\Gamma$
which can be found in \cite[Theorem 1.1]{Li-ci-soluble} and its
proof.

\begin{rem}\label{smallval-con-rem}
{\rm Let  $\Gamma=\Cos(G,H,HgH)$ as in  Construction
\ref{coset-constr}.  Then   $G=\langle H,g\rangle$, $\Gamma$ is
connected and $G$-arc transitive  of valency $p$, $\Aut \Gamma \cong
G\times Z_2$, $|V\Gamma|=|G:H|=2p+2$. Further,
$\diam(\Gamma)=\girth(\Gamma)=3$, so    $\Gamma$ is  not $2$-arc
transitive.

Again, by the proof of \cite[Theorem 1.1]{Li-ci-soluble}, the action
of $\Aut \Gamma$ on $V\Gamma$ has a unique system of  imprimitivity
$\mathcal{B}=\{\Delta_1,\Delta_2,\cdots,\Delta_{p+1}\}$, with
$\Delta_i=\{v_i,v_i'\}$ of size 2, and the kernel of the action of
$\Aut \Gamma$ on $\mathcal{B}$ has order 2. Moreover, $v_i$ is not
adjacent to $v_i'$, and for each $j\neq i$, $v_i$ is adjacent to
exactly one point of $\Delta_j$ and $v_i'$ is adjacent to the other.
Thus, $\Gamma(v_1)\cap \Gamma(v_1')=\emptyset$, $V\Gamma=\{v_1\}\cup
\Gamma(v_1) \cup \{v_1'\}\cup \Gamma(v_1') $, and $\Gamma$ is a
nonbipartite double cover of $K_{p+1}$.}
\end{rem}

\begin{lemma}\label{smallval-p}
The graph $\Gamma=\Cos(G,H,HgH)$ in  Construction \ref{coset-constr}
is  $2$-geodesic transitive but not $2$-arc transitive.
\end{lemma}
{\bf Proof.} Let $A:=\Aut \Gamma$, $v_1\in V\Gamma$ and $u\in
\Gamma(v_1)$. Let $E$ be the kernel of the $A$-action on  $\mathcal
{B}$ and $\overline{A}=A/E$.   Then by the proof of \cite[Theorem
1.1]{Li-ci-soluble}, $E\cong Z_2 \lhd A$, $A=G\times E$,
$\overline{A}\cong G=PSL(2,p)$ and $\overline{A}_{\Delta_1}\cong
A_{v_1}$.  Since $A\cong G\times Z_2$, it follows that
$|A_{v_1}|=\frac{p(p-1)}{2}$, and by Lemma 2.4 of
\cite{Li-ci-soluble}, $A_{v_1}\cong Z_p:Z_{\frac{p-1}{2}}$, which
has a unique permutation action of degree $p$, up to permutational
isomorphism. Since $\Gamma$ is $A$-arc-transitive, $A_{v_1}$ is
transitive on $\Gamma(v_1)$ and hence on $\mathcal{B}\setminus
\{\Delta_1\}$, and therefore also on $\Gamma(v_1')$, all of degree
$p$. Thus the $A_{v_1}$-orbits in $V\Gamma$ are
$\{v_1\},\Gamma(v_1),\Gamma(v_1')$ and $\{v_1'\}$, and it follows
that $\Gamma(v_1')=\Gamma_2(v_1)$. Moreover,  $A_{v_1,u}\cong
Z_{\frac{p-1}{2}}$ has orbit lengths $1,\frac{p-1}{2},\frac{p-1}{2}$
in $\Gamma(v_1)$, and hence has the same orbit lengths in
$\Gamma_2(v_1)$, and also in $\Gamma(u)$ (since $A_{v_1,u}$ is the
point stabiliser of $A_u$ acting on $\Gamma(u)$). Since
$\Gamma(v_1)\cap \Gamma(u)\neq \emptyset$, it follows that the
$A_{v_1,u}$-orbits in $\Gamma(u)$ are $\{v_1\},\Gamma(v_1)\cap
\Gamma(u)$, and $\Gamma_2(v_1)\cap \Gamma(u)$. It follows that
$\Gamma$ is $(A,2)$-geodesic transitive. Since $\girth(\Gamma)=3$,
$\Gamma$ is not 2-arc transitive.  $\Box$

\end{document}